\documentclass[12pt,oneside]{amsart}
\usepackage{amssymb,verbatim,enumerate,amsmath,ifthen,cite}
\usepackage[mathscr]{eucal}
\usepackage[utf8]{inputenc}
\usepackage[T1]{fontenc}
\usepackage[marginparwidth=2.4cm]{geometry}
\newtheorem{theorem}{Theorem}[section]
\newtheorem*{theorem*}{Theorem}
\def\Thm#1#2{\ifthenelse{\equal{#1}{*}}{\begin{theorem*}#2\end{theorem*}}
             {\begin{theorem}\label{T#1}#2\end{theorem}}}
\newtheorem{Atheorem}{Theorem}

\def\THM#1#2{\begin{Atheorem}\label{T#1}#2\end{Atheorem}}
\def\thm#1{Theorem~\ref{T#1}}

\newtheorem{proposition}[theorem]{Proposition}
\newtheorem*{proposition*}{Proposition}
\def\Prp#1#2{\ifthenelse{\equal{#1}{*}}{\begin{proposition*}#2\end{proposition*}}
             {\begin{proposition}\label{P#1}#2\end{proposition}}}
\def\prp#1{Proposition~\ref{P#1}}

\newtheorem{corollary}[theorem]{Corollary}
\newtheorem*{corollary*}{Corollary}
\def\Cor#1#2{\ifthenelse{\equal{#1}{*}}{\begin{corollary*}#2\end{corollary*}}
             {\begin{corollary}\label{C#1}#2\end{corollary}}}
\def\cor#1{Corollary~\ref{C#1}}

\newtheorem{lemma}[theorem]{Lemma}
\newtheorem*{lemma*}{Lemma}
\def\Lem#1#2{\ifthenelse{\equal{#1}{*}}{\begin{lemma*}#2\end{lemma*}}
             {\begin{lemma}\label{L#1}#2\end{lemma}}}
\def\lem#1{Lemma~\ref{L#1}}
\newtheorem{Alemma}{Lemma}

\theoremstyle{definition}
\newtheorem{remark}[theorem]{Remark}
\newtheorem*{remark*}{Remark}
\def\Rem#1#2{\ifthenelse{\equal{#1}{*}}{\begin{remark*}\rm #2\end{remark*}}
             {\begin{remark}\label{R#1}\rm #2\end{remark}}}

\newtheorem{example}[theorem]{Example}
\newtheorem*{example*}{Example}
\def\Exa#1#2{\ifthenelse{\equal{#1}{*}}{\begin{example*}\rm #2\end{example*}}
             {\begin{example}\label{Ex#1}\rm #2\end{example}}}

\numberwithin{equation}{section}
\def\eq#1{{\rm(\ref{#1})}}
\def\Eq#1#2{\ifthenelse{\equal{#1}{*}}
  {\begin{equation*}\begin{aligned}[]#2\end{aligned}\end{equation*}}
  {\begin{equation}\begin{aligned}[]\label{#1}#2\end{aligned}\end{equation}}}

\newcommand{\operator}[1]{\mathop{\vphantom{\sum}\mathchoice
{\vcenter{\hbox{\LARGE $#1$}}}
{\vcenter{\hbox{\large $#1$}}}{#1}{#1}}\displaylimits}

\def\Mst_#1^#2{\operator{\mathscr{M}_{\mbox{\scriptsize$\#$}}\!\!}_{#1}^{#2}\,\,}

\newcommand\R{\mathbb{R}}

\newcommand\N{\mathbb{N}}

\newcommand\Q{\mathbb{Q}}

\newcommand{\floor}[1]{\left\lfloor #1 \right\rfloor}

\DeclareMathOperator{\conv}{conv}

\def\comment#1{}

\newcounter{LT}

\title[Equality and comparison of generalized quasiarithmetic means]{Equality and comparison of generalized quasiarithmetic means}

\author{Zsolt P\'ales}
\address{Institute of Mathematics, University of Debrecen, Pf.\ 400, 4002 Debrecen, Hungary}
\email{pales@science.unideb.hu}

\author{Pawe{\l} Pasteczka}
\address{Institute of Mathematics, University of the National Education Commission, Podcho\-r\k{a}\.{z}ych str 2, 30-084 Krak\'ow, Poland}
\email{pawel.pasteczka@uken.krakow.pl}

\thanks{The first author was supported by the K-134191 NKFIH Grant.}



\begin{document}
\begin{abstract}
The purpose of this paper is to extend the definition of quasiarithmetic means by taking a strictly monotone generating function instead of a strictly monotone and continuous one. We establish the properties of such means and compare them to the analogous properties of standard quasiarithmetic means.
The comparability and equality problems of generalized  quasiarithmetic are also solved. We also provide an example of a mean which, depending on the underlying interval or on the number of variables, could be or could not be represented as a generalized  quasiarithmetic mean.
\end{abstract}
\maketitle

\section{Introduction}

Let $I$ denote a nonempty open subinterval of $\R$ in this paper.
In what follows, for an $n$-tuple $v\in\R^n$, the $i$th coordinate of $v$ will be denoted by $v_i$ ($i\in\{1,\dots,n\}$). Given a strictly increasing continuous function $f:I\to\R$ and $n\in\N$, the \emph{$n$-variable quasiarithmetic mean} $A^{[n]}_{f}:I^n\to I$ is defined by the following formula:
\Eq{*}{
	A^{[n]}_{f}(x):=f^{-1}\bigg(\frac{f(x_1)+\dots+f(x_n)}{n}\bigg) \qquad (x\in I^n).
}
To introduce the notion of the weighted quasiarithmetic mean, denote
\Eq{*}{
  \Lambda_n:=\big\{\lambda\in[0,\infty)^n\,\big|\,\lambda_1+\dots+\lambda_n>0\big\}.
}
Then, \emph{the $n$-variable weighted quasiarithmetic mean} $\widetilde{A}^{[n]}_{f}:I^n\times\Lambda_n\to I$ is defined by 
\Eq{*}{
	\widetilde{A}^{[n]}_{f}(x;\lambda):=f^{-1}\bigg(\frac{\lambda_1f(x_1)+\dots+\lambda_nf(x_n)}{\lambda_1+\dots+\lambda_n}\bigg) \qquad (x\in I^n,\,\lambda\in\Lambda_n).
}
It is well-known that these means possess the mean value property, that is, the following inequalities hold:
\Eq{MV}{
    \min(x)\leq A^{[n]}_{f}&(x)\leq \max(x)\qquad (x\in I^n)
    \qquad\mbox{and}\\
    \min\{x_i:\lambda_i>0\}\leq \widetilde{A}^{[n]}_{f}(&x;\lambda)\leq \max\{x_i:\lambda_i>0\}\qquad (x\in I^n,\,\lambda\in\Lambda_n).
}

The theory of these means was elaborated in several monographs, for instance in \cite{HarLitPol34}. (See also \cite{Bul03}, \cite{BulMitVas88}, \cite{MitPecFin93}.) We recall three important results which motivated this paper. The first result establishes numerous equivalent conditions for the equality of two quasiarithmetic means.

\THM{A}{Let $f,g:I\to\R$ be strictly increasing continuous functions.
Then the following assertions are equivalent:
\begin{enumerate}[(i)]
 \item There exists $n\in\N\setminus\{1\}$ such that, for all $x\in I^n$,
 \Eq{*}{
   A^{[n]}_{f}(x)=A^{[n]}_{g}(x).
 }
 \item For all $n\in\N$ and $x\in I^n$,
 \Eq{*}{
   A^{[n]}_{f}(x)=A^{[n]}_{g}(x).
 }
 \item For all $n\in\N$, $x\in I^n$, and $\lambda\in\Lambda_n$
 \Eq{*}{
   \widetilde{A}^{[n]}_{f}(x;\lambda)=\widetilde{A}^{[n]}_{g}(x;\lambda).
 }
 \item There exists $n\in\N\setminus\{1\}$ such that, for all $x\in I^n$ and $\lambda\in\Lambda_n$
 \Eq{*}{
   \widetilde{A}^{[n]}_{f}(x;\lambda)=\widetilde{A}^{[n]}_{g}(x;\lambda).
 }
 \item There exist two real constants $\alpha>0$ and $\beta$, such that
 \Eq{*}{
   g=\alpha f+\beta.
 }
\end{enumerate}}

The following result offers equivalent conditions for the comparability of two quasiarithmetic means.

\THM{B}{Let $f,g:I\to\R$ be strictly increasing continuous functions. Then the following assertions are equivalent:
\begin{enumerate}[(i)]
 \item There exists $n\in\N\setminus\{1\}$ such that, for all $x\in I^n$,
 \Eq{*}{
   A^{[n]}_{f}(x)\leq A^{[n]}_{g}(x).
 }
 \item For all $n\in\N$ and $x\in I^n$,
 \Eq{*}{
   A^{[n]}_{f}(x)\leq A^{[n]}_{g}(x).
 }
 \item For all $n\in\N$, $x\in I^n$, and $\lambda\in\Lambda_n$
 \Eq{*}{
   \widetilde{A}^{[n]}_{f}(x;\lambda)\leq \widetilde{A}^{[n]}_{g}(x;\lambda).
 }
 \item There exists $n\in\N\setminus\{1\}$ such that, for all $x\in I^n$ and $\lambda\in\Lambda_n$
 \Eq{*}{
   \widetilde{A}^{[n]}_{f}(x;\lambda)\leq \widetilde{A}^{[n]}_{g}(x;\lambda).
 }
 \item For all $x<t<y$ in $I$,
 \Eq{*}{
  \frac{f(y)-f(t)}{f(y)-f(x)}\leq\frac{g(y)-g(t)}{g(y)-g(x)}.
 }
 \item There exists a positive function $p:I\to\R_+$ such that, for all $x,t\in I$
 \Eq{*}{
   f(x)-f(t)\leq p(t)(g(x)-g(t)).
 }
 \item The function $g\circ f^{-1}$ is convex on the interval $f(I)$.
 \item The function $f\circ g^{-1}$ is concave on the interval $g(I)$.
\end{enumerate}}

Finally, we recall the characterization theorem of quasiarithmetic means which is due to Kolmogorov \cite{Kol30}.

\THM{C}{Let $M:=\big(M^{[n]}\big)_{n\in\N}$ be a sequence of functions, where $M^{[n]}:I^n\to I$. Then $M$ is generated by a quasiarithmetic mean, i.e., there exists a strictly increasing continuous function $f:I\to\R$ such that, for all $n\in\N$ and $x\in I^n$,
\Eq{*}{
  M^{[n]}(x)=A^{[n]}_{f}(x)
}
if and only if $M$ possesses the following properties:
\begin{enumerate}[(i)]
 \item $M$ is reflexive, i.e., for all $n\in\N$ and $x\in I^n$ such that $x_1=\dots=x_n$, 
 \Eq{*}{
   M^{[n]}(x)=x_1.
 }
 \item $M$ is continuous and symmetric, i.e., for all $n\in\N$, the function $M^{[n]}$ is continuous and symmetric.
 \item $M$ is strictly increasing, i.e., for all $n\in\N$, for all $i\in\{1,\dots,n\}$ and for all fixed $x_1,\dots,x_{i-1},x_{i+1},\dots,x_n\in I$, the map 
  \Eq{*}{
   x_i\mapsto M^{[n]}(x_1,\dots,x_n) 
  }
  is strictly increasing on $I$.
  \item $M$ is associative, i.e., for all $n\in\N$, for all $k\in\{1,\dots,n-1\}$ and $x\in I^n$, 
 \Eq{*}{
  M^{[n]}(x)=M^{[n]}(y,\dots,y,x_{k+1},\dots,x_n),
 }
 where $y:=M^{[k]}(x_1,\dots,x_k)$.
\end{enumerate}}

We note that there are other known characterizations of quasiarithmetic means due to Aczél \cite{Acz47b,Acz48a,Acz56c}, Acz\'el--Maksa \cite{AczMak96a,AczMak96b}, Burai--Kiss--Szokol \cite{BurKisSzo21}, de Finetti \cite{Def31}, and Nagumo \cite{Nag31}.

We also note that the notion of quasiarithmetic mean has been extended in several ways. One early generalization was introduced by Bajraktarevi\'c \cite{Baj58,Baj63}. Then, Daróczy introduced the concepts of a deviation function and the related deviation mean in \cite{Dar71b,Dar72b} which was further generalized to the concepts of quasideviation function and quasideviation mean by P\'ales in\cite{Pal82a}. In fact, in this paper and later in \cite{Pal89b}, a Kolmogorov-type characterization was established. The notion of semideviation function and the related semideviation mean showed up in the paper \cite{PalPas19a}. Nonsymmetric generalizations of quasiarithmetic means were introduced by Matkowski \cite{Mat10b} and by Gr\"unwald--P\'ales \cite{GruPal20}. A further extension of these concepts and results was obtained by Barczy--P\'ales in \cite{BarPal24a,BarPal24c}.

The purpose of this paper is to extend the definition of quasiarithmetic means by taking a strictly monotone generating function instead of a strictly monotone and continuous one. We establish the properties of such means and then compare them to the analogous properties of standard quasiarithmetic means.
In Section 3, we characterize the comparability of two generalized  quasiarithmetic means, while in Section 4, we deal with their equality problem. In the last section, we provide an example of a mean which depending on the underlying interval or on the number of variables could be or could not be represented as a generalized  quasiarithmetic mean.

\section{Basic properties of generalized quasiarithmetic means}

For a locally bounded function $f:I\to\R$, define $f_-,f_+:I\to\R$ by
\Eq{*}{
  f_-(x):=\lim_{r\to0+}\inf_{u\in(x-r,x+r)\cap I}f(u)
  \quad\mbox{and}\quad
  f_+(x):=\lim_{r\to0+}\sup_{u\in(x-r,x+r)\cap I}f(u)\qquad(x\in I).
}
Then, one can easily see that $f_-$ is the largest lower semicontinuous function below $f$ and $f_+$ is the smallest upper semicontinuous function above $f$, which we call the lower and upper semicontinuous envelopes of $f$, respectively. 

If $f:I\to\R$ is a strictly increasing continuous function, then its range $f(I)$ is also an open interval and $f$ is a homeomorphism between $I$ and $f(I)$, i.e., the inverse of $f$ is continuous, strictly increasing and maps $f(I)$ onto $I$. Furthermore, we have that $f_-=f_+=f$.

In the more general case, when $f:I\to\R$ is a strictly increasing function, then $f:I\to f(I)$ is a bijection and $f$ is also locally bounded. Due to the opennes of the interval $I$, the function $f$ has a left and a right limit at every point of $I$, moreover, we have that
\Eq{*}{
  f_-(x)=\lim_{u\to x-}f(u)
  \qquad\mbox{and}\qquad
  f_+(x):=\lim_{u\to x+}f(u)\qquad(x\in I).
}
One can also easily see that $f_-,f_+:I\to\R$ are strictly increasing functions and, for all elements $u<x<v$ in $I$, we have that
\Eq{*}{
   f_+(u)<f_-(x)\leq f(x)\leq f_+(x)<f_-(v).
}
In the sequel, the set of points in $I$, where $f$ is continuous, will be denoted by $C_f$. It is well-known that $C_f$ is co-countable in $I$, i.e., $I\setminus C_f$ is countable. It is also clear that $x$ belongs to $C_f$ if and only if $f_-(x)=f(x)=f_+(x)$. 
For a strictly decreasing function, similar properties can be established. 

For a subset $S\subseteq\R$, the smallest convex set containing $S$, which is identical to the smallest interval containing it, will be denoted by $\conv(S)$. For our definition of generalized quasiarithmetic means, we shall need the following lemma about the existence and properties of the generalized inverse of strictly increasing (but not necessarily continuous) functions.

\Lem{SMF}{{\rm(\!\!\cite[Lemma 1]{GruPal20})} Let $f:I\to\R$ be a strictly increasing function. Then there exists a uniquely determined increasing function $f^{(-1)}:\conv(f(I))\to I$ which is the left inverse of $f$, i.e., for all $x\in I$,
\Eq{SMF}{
   (f^{(-1)}\circ f)(x)=x.
}
Furthermore, we have the following assertions:
\begin{enumerate}[(i)]
 \item $f^{(-1)}$ is continuous;
 \item For all $y\in f(I)$,
 \Eq{*}{
  (f\circ f^{(-1)})(y)=y;
 }
 \item For all $y\in\conv(f(I))$,
 \Eq{*}{
  (f_-\circ f^{(-1)})(y)\leq y \leq (f_+\circ f^{(-1)})(y);
 }
 \item For all $x\in I$,
 \Eq{*}{
   (f^{(-1)}\circ f_-)(x)=(f^{(-1)}\circ f_+)(x)=x,
 }
 equivalently,
 \Eq{*}{
   (f_-)^{(-1)}=(f_+)^{(-1)}=f^{(-1)}.
 }
\end{enumerate}
}

We should note that assertion (iv) was not explicitly stated in \cite[Lemma 1]{GruPal20}, but it easily follows from \eq{SMF} by taking left and right limits at $x$ and using assertion (i).

The function $f^{(-1)}$ described in the above lemma is called the \emph{generalized inverse of the strictly increasing function $f:I\to\R$}. It is clear from assertions (ii) and (iii) that the restriction of $f^{(-1)}$ to $f(I)$ equals the inverse of $f$ in the standard sense. Therefore, $f^{(-1)}$ is the continuous and increasing extension of the inverse of $f$ to the smallest interval containing the range of $f$.

As an easy consequence of this lemma, we have the following useful statement which will not require a proof.

\Lem{ineq}{Let $f:I\to\R$ be a strictly increasing function and let $x\in I$ and $u\in\conv(f(I))$. Then we have the following equivalences:
\begin{enumerate}[(i)]
 \item $f^{(-1)}(u)=x$ holds if and only if $u\in[f_-(x),f_+(x)]$.
 \item $f^{(-1)}(u)<x$ holds if and only if $u<f_-(x)$.
 \item $f^{(-1)}(u)\leq x$ holds if and only if $u\leq f_+(x)$.
 \item $x<f^{(-1)}(u)$ holds if and only if $f_+(x)<u$.
 \item $x\leq f^{(-1)}(u)$ holds if and only if $f_-(x)\leq u$.
\end{enumerate}
}

Given a strictly increasing function $f:I\to\R$ and $n\in\N$, the \emph{$n$-variable generalized quasiarithmetic mean} $A^{[n]}_{f}:I^n\to I$ and the \emph{$n$-variable weighted generalized quasiarithmetic mean} $\widetilde{A}^{[n]}_{f}:I^n\times\Lambda_n\to I$ are defined by the following formulas:
\Eq{QAM}{
	A^{[n]}_{f}(x):=f^{(-1)}\bigg(\frac{f(x_1)+\dots+f(x_n)}{n}\bigg) \qquad (x\in I^n),
}
\Eq{WQAM}{
	\widetilde{A}^{[n]}_{f}(x;\lambda):=f^{(-1)}\bigg(\frac{\lambda_1f(x_1)+\dots+\lambda_nf(x_n)}{\lambda_1+\dots+\lambda_n}\bigg) \qquad (x\in I^n,\,\lambda\in\Lambda_n).
}
Clearly, for any $x\in I^n$, we have that $A^{[n]}_{f}(x)=\widetilde{A}^{[n]}_{f}(x;\lambda)$ if $\lambda_1=\dots=\lambda_n$.
We also mention that, according to Example 4.5(3) of the paper \cite{PalPas19a}, the mean $\widetilde{A}^{[n]}_{f}$ turns out to be a weighted semideviation mean generated by the semideviation function  $D_f(x,t):=f(x)-f(t)$ \ ($x,t\in I$). On the other hand, these means are the particular cases of the generalized Bajraktarevi\'c means which were introduced in the paper \cite{GruPal20}. 
The following result was established in this latter paper.

\Prp{QAM}{{\rm(\!\!\cite[Theorem 2]{GruPal20})} Let $f:I\to\R$ be strictly increasing and let $n\in\N$. Then the functions $A^{[n]}_{f}:I^n\to I$ and $\widetilde{A}^{[n]}_{f}:I^n\times\Lambda_n\to I$ given by \eq{QAM} and \eq{WQAM} are well-defined and they possess the mean value property, that is, \eq{MV} holds.}

We should note that the second inequality was not explicitly stated in \cite{GruPal20}, but it can be proved similarly as the first one.
Further properties of these means are obtained in the following result.

\Prp{QAM2}{Let $f:I\to\R$ be strictly increasing and let $n\in\N$. Then we have the following assertions:
\begin{enumerate}[(i)]
 \item The means $A^{[n]}_{f}$ and $\widetilde{A}^{[n]}_{f}$ are symmetric functions, i.e., for all permutation $\pi$ of the index set $\{1,\dots,n\}$, for all $x\in I^n$ and $\lambda\in\Lambda_n$, we have that
 \Eq{*}{
   A^{[n]}_{f}(x_1,\dots,x_n)&=A^{[n]}_{f}(x_{\pi(1)},\dots,x_{\pi(n)}) \qquad\mbox{and}\\
   \widetilde{A}^{[n]}_{f}(x_1,\dots,x_n;\lambda_1,\dots,\lambda_n)
   &=\widetilde{A}^{[n]}_{f}(x_{\pi(1)},\dots,x_{\pi(n)};\lambda_{\pi(1)},\dots,\lambda_{\pi(n)}).
 }
 \item The means $A^{[n]}_{f}$ and $\widetilde{A}^{[n]}_{f}$ are increasing functions of their variables, i.e., for all $i\in\{1,\dots,n\}$, for all fixed $x_1,\dots,x_{i-1},x_{i+1},\dots,x_n\in I$ and $\lambda\in\Lambda_n$, the maps 
  \Eq{*}{
   x_i\mapsto A^{[n]}_{f}(x_1,\dots,x_n) \qquad\mbox{and}\qquad
   x_i\mapsto\widetilde{A}^{[n]}_{f}(x_1,\dots,x_n;\lambda_1,\dots,\lambda_n)
  }
  are increasing on $I$.
  \item The mean $\widetilde{A}^{[n]}_{f}$ is a continuous function of its weights, i.e., for all fixed $x\in I^n$, the map
  \Eq{*}{
    \lambda\mapsto \widetilde{A}^{[n]}_{f}(x;\lambda)
  }
  is continuous over $\Lambda_n$.
\end{enumerate}}

\begin{proof}
Assertion (i) is direct consequence of the definition of the means $A^{[n]}_{f}$ and $\widetilde{A}^{[n]}_{f}$. Using that $f$ and $f^{(-1)}$ are increasing functions, assertion (ii) follows easily.
Applying the continuity of $f^{(-1)}$, we can also see that assertion (iii) is also valid. 
\end{proof}

The following lemma will be useful in the sequel.

\Lem{SC}{Let $f:I\to\R$ be strictly increasing, $x\in I$, and $n\in\N\setminus\{1\}$. Then we have the following assertions:
\begin{enumerate}[(i)]
 \item $f$ is lower semicontinuous at $x$ if and only if, for all $y\in I$ with $y<x$,
\Eq{i}{
  A^{[n]}_{f}(\underbrace{x,\dots,x}_{n-1 \text{ times}},y)<x.
}
\item $f$ is upper semicontinuous at $x$ if and only if, for all $y\in I$ with $y>x$,
\Eq{ii}{
  A^{[n]}_{f}(\underbrace{x,\dots,x}_{n-1 \text{ times}},y)>x.
}
\item $f$ is continuous at $x$ if and only if, for all $y\in I$ with $y\neq x$
\Eq{iii}{
  A^{[n]}_{f}(\underbrace{x,\dots,x}_{n-1 \text{ times}},y)\neq x.
}
\end{enumerate}}

\begin{proof}
Assume that $f$ is lower semicontinuous at $x$. Then, for all $y<x$ in $I$, we have that $f(y)<f_-(x)=f(x)$. Therefore, 
\Eq{*}{
  \frac{(n-1)f(x)+f(y)}{n}<f(x)=f_-(x),
}
which, in view of assertion (ii) of \lem{ineq}, implies the inequality \eq{i}.

Conversely, assume that the inequality \eq{i} holds for all $y<x$ in $I$. Then, again by assertion (ii) of \lem{ineq}, we get that
\Eq{*}{
  \frac{(n-1)f(x)+f(y)}{n}<f_-(x).
}
Upon taking the limit as $y$ tends to $x$ from the left, it follows that
\Eq{*}{
  \frac{(n-1)f(x)+f_-(x)}{n}\leq f_-(x).
}
After rearranging this inequality, we conclude that $f(x)\leq f_-(x)$, which shows that $f$ is lower semicontinuous at the point $x$.

The proof of assertion (ii) is analogous. Assertion (iii) is an easy consequence of assertions (i) and (ii).
\end{proof}

Now we show several equivalent conditions equivalent to the continuity of the generating function $f$. Let us emphasize that conditions (iv), (vi), and (viii) of the theorem below correspond to Kolmogorov's conditions (ii), (iii), and (iv), respectively (see \thm{C}). We also mention that the proof uses \lem{SC} and also \prp{QAM3}, \prp{SC}, \cor{SC} which will be established in Section 3.

\Thm{QAM=AM}{
Let $f:I\to\R$ be a strictly increasing function. Then the following conditions are equivalent:
\begin{enumerate}[\rm (i)]
  \item The function $f$ is continuous on $I$. 
  \item There exists $n\in\N\setminus\{1\}$ and a continuous, strictly increasing function $g \colon I \to \R$ such that $A_{f}^{[n]}=A^{[n]}_{g}$ holds on $I^n$ (i.e., $A_{f}^{[n]}$ is a standard quasiarithmetic mean on $I^n$).
 \item There exists $n\in\N\setminus\{1\}$ and a continuous, strictly increasing function $g \colon I \to \R$ such that either $A^{[n]}_{f}\le A^{[n]}_{g}$ holds on $I^n$ or $A^{[n]}_{g}\le A^{[n]}_{f}$ holds on $I^n$. 
 \item There exists $n\in\N\setminus\{1\}$ such that $A^{[n]}_{f}$ is continuous.
  \item There exists $n\in\N\setminus\{1\}$ such that $A^{[n]}_{f_+}\leq A^{[n]}_{f_-}$ holds on $I^n$.
  \item There exists $n\in\N\setminus\{1\}$ such that  $A^{[n]}_{f}$ is strictly increasing. 
 \item There exists $n\in\N\setminus\{1\}$ such that  $A^{[n]}_{f}$ is strict, that is 
 \Eq{*}{
 \min(x)<A^{[n]}_{f}(x)<\max(x)\qquad \text{ for every non-diagonal vector }x\in I^n.
 }
  \item $(A_{f}^{[n]})_{n \in \N}$ is associative, that is, the sequence of means $M^{[n]}:=A_f^{[n]}$, $n \in \N$, satisfies property (iv) in \thm{C}. 
\end{enumerate}
Moreover, if any of the above conditions holds, then {\rm (iii)--(vii)} are valid for all $n \in\N$.
}
\begin{proof}
First note that if $f$ is continuous, then all the other listed conditions are satisfied by the characterization theorem of Kolmogorov. 

Obviously, (ii) implies (iii). On the other hand, if (iii) holds, then, according to statement (iii) of \prp{SC}, we get that $C_f=C_g=I$, and hence $f$ has to be continuous, i.e., condition (i) follows. Thus, we have proved the equivalence of conditions (i), (ii), and (iii). 

If (iv) holds, then, by \prp{QAM3}, we have that $A^{[n]}_{f_+}=A^{[n]}_{f_-}=A^{[n]}_{f}$, therefore, condition (v) is satisfied. On the other hand, if (v) is valid, then, according to \cor{SC}, it follows that $f$ is continuous. This shows the equivalence of conditions (i), (iv), and (v).

It is easy to see that (vi) implies (vii). If (vii) is valid, i.e., 
the mean $A^{[n]}_{f}$ is strict, then, according to assertion (iii) of \lem{SC}, it follows that $f$ is continuous on $I$. Thus, we have verified the equivalence of conditions (i), (vi), and (vii).

Finally, we show that condition (viii) implies condition (vii).
To the contrary, assume that (viii) holds true but (vii) is not valid. Then there exists $n \in \N \setminus\{1\}$ and a non-diagonal vector $x \in I^n$ such that either 
$A_f^{[n]}(x)=\max(x)$ or $A_f^{[n]}(x)=\min(x)$. 
Without loss of generality, we can assume that $A_f^{[n]}(x)=\max(x)$ holds and set $a:=\min(x)$,  $b:=\max(x)$. Finally, for all $m,k \in\N$, we define the vector
\Eq{*}{
w^{(m,k)}:=(\underbrace{a,\dots,a}_{m\text{ times}},\underbrace{b,\dots,b}_{k\text{ times}}).
}
Since $A_f^{[n]}$ is monotone, $A_f^{[n]}(x)=b$ yields $A_f^{[n]}(w^{(1,n-1)})=b$. Now, we are going to show by induction on $m$ that $A_f^{[m+n-1]}(w^{(m,n-1)})=b$ holds for all $m\in\N$. We have seen this property for $m=1$. Assume now that it holds for some $m\in\N$. Then, by using the associativity of $A_f$ twice, we get 
\Eq{*}{
A_f^{[n+m]}(w^{(m+1,n-1)})&=A_f^{[n+m]}(a,w^{(m,n-1)})\\
&=A_f^{[n+m]}\Big(a,A_f^{[n+m-1]}\big(w^{(m,n-1)}\big),\dots,A_f^{[n+m-1]}\big(w^{(m,n-1)}\big)\Big)\\
&=A_f^{[n+m]}\big(w^{(1,n+m-1)}\big)
=A_f^{[n+m]}\big(w^{(1,n-1)},w^{(0,m)}\big)\\
&=A_f^{[n+m]}\Big(A_f^{[n]}(w^{(1,n-1)}),\dots,A_f^{[n]}(w^{(1,n-1)}),w^{(0,m)}\Big)\\
&=A_f^{[n+m]}\big(w^{(0,n)},w^{(0,m)}\big)=b,
}
which completes the induction. 
Therefore we have $A_f^{[m+n-1]}(w^{(m,n-1)})=b$ for all $m \ge 1$. On the other hand, due to the continuity of $f^{(-1)}$, we get
\Eq{*}{
\lim_{m \to \infty} A_f^{[m+n-1]}(w^{(m,n-1)})
&=\lim_{m \to \infty}f^{(-1)} \Big( \frac{mf(a)+(n-1)f(b)}{m+n-1}\Big)\\
&=f^{(-1)} \Big( \lim_{m \to \infty}\big(\tfrac{m}{m+n-1}f(a)+\tfrac{n-1}{m+n-1}f(b)\big)\Big)\\
&=f^{(-1)} \big(  f(a)\big)=a,
}
which leads to the contradiction $b=a$.
This completes the proof of the final implication.
\end{proof}

The content of the above theorem can be summarized as follows: If a generalized quasiarithmetic mean is not equal (or not comparable) to a standard quasiarithmetic mean, then it is not continuous, not strictly increasing (but increasing), not a strict mean and it is also not associative. 

In the next proposition, we establish the properties of the generalized quasiarithmetic mean generated by the lower and upper semicontinuous envelopes of $f$.

\Prp{QAM3}{Let $f:I\to\R$ be strictly increasing and let $n\in\N$. Then we have the following assertions:
\begin{enumerate}[(i)]
  \item For $x\in I^n$ and $\lambda\in\Lambda_n$, we have that
 \Eq{*}{
   \lim_{u\to x-}A^{[n]}_{f}(u)
   =A^{[n]}_{f_-}(x)
   \leq A^{[n]}_{f}(x)
   \leq A^{[n]}_{f_+}(x)
   =\lim_{u\to x+}A^{[n]}_{f}(u)\qquad\mbox{and}\\
  \lim_{u\to x-}\widetilde{A}^{[n]}_{f}(u;\lambda)
   =\widetilde{A}^{[n]}_{f_-}(x;\lambda)
   \leq \widetilde{A}^{[n]}_{f}(x;\lambda)
   \leq \widetilde{A}^{[n]}_{f_+}(x;\lambda)
   =\lim_{u\to x+}\widetilde{A}^{[n]}_{f}(u;\lambda).
 }
 \item The means $A^{[n]}_{f_-}$ and $\widetilde{A}^{[n]}_{f_-}$ are lower semicontinuous and the means $A^{[n]}_{f_+}$ and $\widetilde{A}^{[n]}_{f_+}$ are upper semicontinuous on their domains. 
 \item The means $A^{[n]}_{f_-}$ and $\widetilde{A}^{[n]}_{f_-}$ are subassociative, i.e., for all $k\in\{1,\dots,n-1\}$ and $x\in I^n$, $\lambda\in\Lambda_n$ with $\lambda_1+\dots+\lambda_k>0$, 
 \Eq{*}{
   A^{[n]}_{f_-}(x)&\geq A^{[n]}_{f_-}(y,\dots,y,x_{k+1},\dots,x_n)\qquad\mbox{and}\\ 
   \widetilde{A}^{[n]}_{f_-}(x;\lambda)&\geq \widetilde{A}^{[n]}_{f_-}(z,\dots,z,x_{k+1},\dots,x_n;\lambda_1,\dots,\lambda_n)
 }
 where $y:=A^{[k]}_{f_-}(x_1,\dots,x_k)$ and $z:=\widetilde{A}^{[k]}_{f_-}(x_1,\dots,x_k;\lambda_1,\dots,\lambda_k)$.
\item The means $A^{[n]}_{f_+}$ and $\widetilde{A}^{[n]}_{f_+}$ are superassociative, i.e., for all $k\in\{1,\dots,n-1\}$ and $x\in I^n$, $\lambda\in\Lambda_n$ with $\lambda_1+\dots+\lambda_k>0$, 
 \Eq{*}{
   A^{[n]}_{f_+}(x)&\leq A^{[n]}_{f_+}(y,\dots,y,x_{k+1},\dots,x_n)\qquad\mbox{and}\\ 
   \widetilde{A}^{[n]}_{f_+}(x;\lambda)&\leq \widetilde{A}^{[n]}_{f_+}(z,\dots,z,x_{k+1},\dots,x_n;\lambda_1,\dots,\lambda_n)
 }
 where $y:=A^{[k]}_{f_+}(x_1,\dots,x_k)$ and $z:=\widetilde{A}^{[k]}_{f_+}(x_1,\dots,x_k;\lambda_1,\dots,\lambda_k)$.
\end{enumerate}}

\begin{proof}
To prove assertion (i) for the weighted mean, let $x\in I^n$ and $\lambda\in\Lambda_n$ and denote $\lambda_0:=\lambda_1+\dots+\lambda_n$.
Then, using assertions (i) and (iv) of \lem{SMF}, we get that
\Eq{*}{
  \lim_{u\to x-}\widetilde{A}^{[n]}_{f}(u;\lambda)
  &=\lim_{u\to x-}f^{(-1)}\bigg(\frac1{\lambda_0}\sum_{i=1}^n\lambda_if(u_i)\bigg)
  =f^{(-1)}\bigg(\frac1{\lambda_0}\sum_{i=1}^n \lambda_i\lim_{u_i\to x_i-}f(u_i)\bigg)\\
  &=f^{(-1)}\bigg(\frac1{\lambda_0}\sum_{i=1}^n \lambda_if_-(x_i)\bigg)
  =\widetilde{A}^{[n]}_{f_-}(x;\lambda),
}
which proves the first equality in the second chain of inequalities-equalities of assertion (i). The verification of the last equality is analogous. The two inequalities in this chain follow directly from the inequalities $f_-\leq f\leq f_+$ and the increasingness of $f^{(-1)}$.
By taking $\lambda_1=\dots=\lambda_n=1$ in the second chain, we obtain the first one. 

By the lower semicontinuity of the function $f_-$, we have that the maps
\Eq{*}{
  x\mapsto \frac1n\sum_{i=1}^nf_-(x_i)\qquad\mbox{and}\qquad
  (x;\lambda)\mapsto \frac1{\lambda_1+\dots+\lambda_n}\sum_{i=1}^n\lambda_if_-(x_i)
}
are lower semicontinuous on $I^n$ and on $I^n\times\Lambda_n$, respectively. Therefore, the composition of these maps with $f^{(-1)}$ is also lower semicontinuous on the respective domains. This proves that assertion (ii) is valid.

To prove the subassociativity of the weighted mean $\widetilde{A}^{[n]}_{f_-}$, let $k\in\{1,\dots,n-1\}$ and $x\in I^n$, $\lambda\in\Lambda_n$ with $\lambda_1+\dots+\lambda_k>0$ and let $z:=\widetilde{A}^{[k]}_{f_-}(x_1,\dots,x_k;\lambda_1,\dots,\lambda_k)$.
Then, using assertion (iii) of \lem{SMF}, we obtain
\Eq{*}{
  f_-(z)
  =f_-\circ f^{(-1)}\bigg(\frac1{\lambda_1+\dots+\lambda_k}\sum_{i=1}^k\lambda_if_-(x_i)\bigg)
  \leq \frac1{\lambda_1+\dots+\lambda_k}\sum_{i=1}^k\lambda_if_-(x_i).
}
Therefore,
\Eq{*}{
  (\lambda_1+\dots+\lambda_k)f_-(z)+\sum_{i=k+1}^n\lambda_if_-(x_i)\leq \sum_{i=1}^n\lambda_if_-(x_i).
}
Substituting the left and right hand sides into $f^{(-1)}$, it follows that the second inequality of assertion (iii) is valid, which proves the subassociativity of $\widetilde{A}^{[n]}_{f_-}$. By taking $\lambda_1=\dots=\lambda_n=1$, the subassociativity of $A^{[n]}_{f_-}$ follows. Thus the proof of assertion (iii) is complete.

The proof of the superassociativity of $A^{[n]}_{f_+}$ and $\widetilde{A}^{[n]}_{f_+}$ is completely analogous and therefore, it is omitted.
\end{proof}

The last result of this section shows that an $n$-variable generalized quasiarithmetic mean completely determines the corresponding $m$-variable generalized quasiarithmetic mean whenever $m<n$.

\Thm{Emn}{
Let $f:I\to\R$ be a strictly increasing function and let $n,m\in\N$ with $m<n$. Then, for all $x \in I^m$,
\Eq{*}{
A^{[m]}_{f}(x_1,\dots,x_m)
&=\inf\{z  \in I \colon A^{[n]}_{f}(x_1,\dots,x_m,z,\dots,z)<z\}\\
&=\sup\{z  \in I \colon A^{[n]}_{f}(x_1,\dots,x_m,z,\dots,z)>z\}.
}
}

\begin{proof} 
We are going to show only the first equality as the second one is analogous. Let $x \in I^m$ be fixed and define the set $M(x)$ by
\Eq{*}{
  M(x):=\{z  \in I \colon A^{[n]}_{f}(x_1,\dots,x_m,z,\dots,z)<z\}.
}
We need to show the equality
\Eq{*}{
  A^{[m]}_{f}(x)=\inf M(x).
}

Let $z\in M(x)$ be arbitrary. Then, by \lem{ineq}, we have 
\Eq{*}{
\frac{f(x_1)+\dots+f(x_m)+(n-m)f(z)}n<f_-(z).
}
Therefore, using also the inequality $f_-(z)\leq f(z)$, we get
\Eq{*}{
\frac{f(x_1)+\dots+f(x_m)}m<\frac{nf_-(z)-(n-m)f(z)}m \le f_-(z),
}
which, again by \lem{ineq}, implies that $A^{[m]}_{f}(x_1,\dots,x_m)<z$.
This proves the inequality $A^{[m]}_{f}(x)\leq\inf M(x)$.

To show the reverse inequality, take $z \in C_f \cap (A^{[m]}_{f}(x_1,\dots,x_m),\infty \big)$. 
Then, we have
\Eq{*}{
f(z) > \frac{f(x_1)+\dots+f(x_m)}m,
}
which yields
\Eq{*}{
f(z) > \frac{f(x_1)+\dots+f(x_m)+(n-m)f(z)}n.
}
Whence, in view of \lem{ineq}, it follows that $A^{[n]}_{f}(x_1,\dots,x_m,z,\dots,z) < z$, and thus  $z \in M(x)$. Therefore $M(x)\supseteq C_f \cap (A^{[m]}_{f}(x_1,\dots,x_m),\infty \big)$ and whence, due to the density of $C_f$, we obtain
\Eq{*}{
\inf M(x) \le \inf C_f \cap (A^{[m]}_{f}(x_1,\dots,x_m),\infty \big) = A^{[m]}_{f}(x_1,\dots,x_m),
}
which completes the proof of the equality $A^{[m]}_{f}(x)=\inf M(x)$.
\end{proof}

\section{Comparison of generalized quasiarithmetic means}

In the result below we establish that, whenever the comparison inequality holds between two $n$-variable generalized quasiarithmetic means and $m<n$, then it holds between the corresponding $m$-variable generalized quasiarithmetic means, as well.

\Thm{M1}{Let $f,g:I\to\R$ be strictly increasing functions and let $n,m\in\N$ with $m<n$. If, for all $x\in I^n$, the inequality
\Eq{nv}{
  A_f^{[n]}(x)\leq A_g^{[n]}(x)
}
holds, then, for all $y\in I^m$,
\Eq{mv}{
  A_f^{[m]}(y)\leq A_g^{[m]}(y).
}}

\begin{proof}
Fix $y=(y_1,\dots,y_m)\in I^m$ arbitrarily. Then, by \eq{nv}, for all $z \in I$, we have
\Eq{*}{
A^{[n]}_{f}(y_1,\dots,y_m,z,\dots,z) \le A^{[n]}_{g}(y_1,\dots,y_m,z,\dots,z).
}
Whence it follows that
\Eq{*}{
\{z  \in I \colon A^{[n]}_{f}(y_1,\dots,y_m,z,\dots,z)<z\} \supseteq \{z  \in I \colon A^{[n]}_{g}(y_1,\dots,y_m,z,\dots,z)<z\}.
}
Upon taking infimum side-by-side, in view of \thm{Emn}, we obtain \eq{mv}.
\end{proof}

In the subsequent result, we establish a necessary condition for the validity of the inequality \eq{nv}. Recall that $\floor{u}$ denotes the lower integer part of $u\in\R$.

\Prp{nv}{Let $f,g:I\to\R$ be strictly increasing functions, $n\in\N\setminus\{1\}$, and assume that \eq{nv} holds for all $x\in I^n$. Then, for all $m\in\{1,\dots,n\}$, $x,y\in I$ with $x<y$ and $t\in(x,y)\cap C_f\cap C_g$,
\Eq{nvnc}{
 \floor{m\frac{f(y)-f(t)}{f(y)-f(x)}}
 \leq \floor{m\frac{g(y)-g(t)}{g(y)-g(x)}}.
}}

\begin{proof} Let $m\in\{1,\dots,n\}$. 
Let $x,y\in I$ with $x<y$ and $t\in(x,y)\cap C_f\cap C_g$. Contrary to the statement, suppose that \eq{nvnc} is not valid, i.e.,
\Eq{*}{
  \floor{m\frac{g(y)-g(t)}{g(y)-g(x)}}
  <\floor{m\frac{f(y)-f(t)}{f(y)-f(x)}}=:k.
}
Then, obviously, $k\in\{1,\dots,m-1\}$. This inequality implies that
\Eq{*}{
  m\frac{g(y)-g(t)}{g(y)-g(x)}< k
  \leq m\frac{f(y)-f(t)}{f(y)-f(x)}.
}
Rearranging these inequalities, it follows that
\Eq{*}{
  \frac{kg(x)+(m-k)g(y)}{m}<g(t) \qquad\mbox{and}\qquad
  f(t)\leq \frac{kf(x)+(m-k)f(y)}{m}.
}
Using that $t\in C_f\cap C_g$ and applying \lem{ineq} to these inequalities, we conclude that
\Eq{*}{
  A^{[m]}_{g}(\underbrace{x,\dots,x}_{k\text{ times}},\underbrace{y,\dots,y}_{m-k\text{ times}})
  <t\leq A^{[m]}_{f}(\underbrace{x,\dots,x}_{k\text{ times}},\underbrace{y,\dots,y}_{m-k\text{ times}}),
}
which shows that the inequality \eq{mv} is not valid for $\big(\underbrace{x,\dots,x}_{k\text{ times}},\underbrace{y,\dots,y}_{m-k\text{ times}}\big)\in I^m$. Therefore, according to \thm{M1}, the inequality \eq{nv} is also not valid for some vector in $I^n$. This contradiction proves that \eq{nvnc} must be valid on the indicated domain.
\end{proof}

Observe that, for any $u,v\in\R$, the inequality $\floor{u}\leq\floor{v}$ is equivalent to $\floor{u}\leq v$, therefore the necessary condition \eq{nvnc} for the validity of the comparison inequality \eq{nv} can be rewritten as
\Eq{*}{
  \max_{m\in\{1,\dots,n\}}\frac1m\floor{m\frac{f(y)-f(t)}{f(y)-f(x)}}
  \leq \frac{g(y)-g(t)}{g(y)-g(x)}
  \quad(x,y\in I,\, x<y,\, t\in(x,y)\cap C_f\cap C_g).
}

The next result shows that if two generalized quasiarithmetic means are comparable in the weakest sense (i.e., the corresponding $n$-variable means are comparable), then their generating functions have the same discontinuity points.

\Prp{SC}{Let $n\in\N\setminus\{1\}$, let $f,g:I\to\R$ be strictly increasing functions, and let $t\in I$. Assume that \eq{nv} holds for all $x\in I^n$. Then, we have the following assertions:
\begin{enumerate}[(i)]
 \item If $g$ is lower semicontinuous at $t$, then so is $f$.
 \item If $f$ is upper semicontinuous at $t$, then so is $g$.
 \item The function $f$ is continuous at $t$ if and only if $g$ is continuous at $t$.
\end{enumerate}}

\begin{proof} Assume that $f$ is not lower semicontinuous at $t$. Then, according to statement (i) of \lem{SC}, there exists a value $y\in I$ with $y<t$ such that
\Eq{*}{
  t\leq A_f^{[n]}(t,\dots,t,y)
}
Using the inequality \eq{nv} for the vector $x=(t,\dots,t,y)\in I^n$, we get
that
\Eq{*}{
  t\leq A_g^{[n]}(t,\dots,t,y).
}
Therefore, again by statement (i) of \lem{SC}, it follows that $g$ is not lower semicontinuous at $t$ as well. This proves assertion (i).

The similar proof of assertion (ii) is based on statement (ii) of \lem{SC}.

Upon taking right limit in \eq{nv}, for all $x\in I^n$, we get that
\Eq{2v+}{
  A_{f_+}^{[n]}(x)=\lim_{u\to x+}A_{f}^{[n]}(u)
  \leq \lim_{u\to x+}A_{g}^{[n]}(u)=A_{g_+}^{[n]}(x).
}

To prove assertion (iii), we show first that $C_g\subseteq C_f$. Let $t\in C_g$. Then $g_+$ is lower semicontinuous at $t$. In view of the inequality \eq{2v+} and assertion (i) (applied for $(f_+,g_+)$ instead of $(f,g)$), it follows that $f_+$ is also lower semicontinuous at $t$ and hence, $f$ is continuous at $t$. This proves that $C_g\subseteq C_f$.

The proof of the reverse inclusion $C_f\subseteq C_g$ is based on the inequality 
\Eq{*}{
  A_{f_-}^{[n]}(x)=\lim_{u\to x-}A_{f}^{[n]}(u)
  \leq \lim_{u\to x-}A_{g}^{[n]}(u)=A_{g_-}^{[n]}(x)
}
and assertion (ii), otherwise it is completely similar to the proof of $C_g\subseteq C_f$.  
\end{proof}

\Cor{SC}{Let $n\in\N\setminus\{1\}$, let $f,g:I\to\R$ be strictly increasing functions. Assume that 
\Eq{+-}{
  A_{f_+}^{[n]}(x)\leq A_{g_-}^{[n]}(x)
}
holds for all $x\in I^n$. Then, $f$ and $g$ are continuous on $I$ and the means in the above inequality are quasiarithmetic in the standard sense.}

\begin{proof} We apply the previous proposition to the pair of functions $(f_+,g_-)$. Since $g_-$ is lower semicontinuous everywhere in $I$, the inequality \eq{+-} and assertion (i) of this proposition imply that $f_+$ has to be lower semicontinuous as well. Therefore, $f$ is continuous everywhere. 

On the other hand, $f_+$ is upper semicontinuous everywhere, thus the inequality \eq{+-} and assertion (ii) of \prp{SC} yield that $g_-$ is upper semicontinuous and hence $g$ is continuous.
\end{proof}

The following theorem is our main result about the comparison of generalized quasiarithmetic means. It is parallel to the classical comparison theorem stated in \thm{B}.

\Thm{M2}{Let $f,g:I\to\R$ be strictly increasing functions. Then the following assertions are equivalent.
\begin{enumerate}[(i)]
 \item For all $n\in\N$ and $x\in I^n$,
 \Eq{*}{
   A^{[n]}_{f}(x)\leq A^{[n]}_{g}(x).
 }
 \item For all $n\in\N$, $x\in I^n$, and $\lambda\in\Lambda_n$,
 \Eq{*}{
   \widetilde{A}^{[n]}_{f}(x;\lambda)\leq \widetilde{A}^{[n]}_{g}(x;\lambda).
 }
 \item There exists $n\in\N\setminus\{1\}$ such that, for all $x\in I^n$ and $\lambda\in\Lambda_n$,
 \Eq{*}{
   \widetilde{A}^{[n]}_{f}(x;\lambda)\leq \widetilde{A}^{[n]}_{g}(x;\lambda).
 }
 \item For all $x,y\in I$ with $x<y$ and $t\in (x,y)\cap C_f\cap C_g$,
 \Eq{*}{
  \frac{f(y)-f(t)}{f(y)-f(x)}\leq\frac{g(y)-g(t)}{g(y)-g(x)}.
 }
 \item There exists a positive function $p:(C_f\cap C_g)\to\R_+$ such that, for all $x\in I$ and $t\in C_f\cap C_g$,
 \Eq{*}{
   f(x)-f(t)\leq p(t)(g(x)-g(t)).
 }
 \item There exists an increasing convex homeomorphism $\varphi:\conv(f(I))\to\conv(g(I))$ such that
 \Eq{*}{
   \varphi(f(x))\leq g(x)\qquad(x\in I)
   \qquad\mbox{and}\qquad
   \varphi(f(t))=g(t)\qquad(t\in C_f\cap C_g).
 }
\item There exists an increasing concave homeomorphism $\psi:\conv(g(I))\to\conv(f(I))$ such that
 \Eq{*}{
   f(x)\leq \psi(g(x))\qquad(x\in I)
   \qquad\mbox{and}\qquad
   f(t)=\psi(g(t))\qquad(t\in C_f\cap C_g).
 }
\end{enumerate}}

\begin{proof} The proof of the implication (i)$\Rightarrow$(ii) is standard. Assuming assertion (i), first one should prove that the inequality in (ii) holds for all weight vector $\lambda$ with rational coordinates, i.e., for all $\lambda\in\Lambda_n\cap\Q^n$. Then, using the continuity of the weighted generalized quasiarithmetic means in their weights (cf.\ assertion (iii) of \prp{QAM2}), this inequality can be verified for all $\lambda\in\Lambda_n$.

The implication (ii)$\Rightarrow$(iii) is obvious.

To verify the implication (iii)$\Rightarrow$(iv), assume that assertion (iii) is valid for some $n\geq 2$. Let $x,y\in I$ with $x<y$ and $\lambda\in (0,1)$. Then applying the inequality in (iii) to the vector
$(x,y,\dots,y)\in I^n$ and to the weight vector $(\lambda,1-\lambda,0,\dots,0)\in\Lambda_n$, we get that
\Eq{2v}{
  \widetilde{A}^{[2]}_{f}(x,y;\lambda,1-\lambda)
  \leq \widetilde{A}^{[2]}_{g}(x,y;\lambda,1-\lambda).
}
Let $t\in(x,y)\cap C_f\cap C_g$ be arbitrary. We prove the inequality asserted in (iv) by contradiction. Both sides of this inequality belong to the interval $[0,1]$. If this inequality is not true, then there exist $\lambda\in(0,1)$ such that
\Eq{*}{
  \frac{g(y)-g(t)}{g(y)-g(x)}<\lambda<\frac{f(y)-f(t)}{f(y)-f(x)}.
}
The left and the right hand side inequalities imply that
\Eq{*}{
  \lambda g(x)+(1-\lambda)g(y)<g(t) \qquad\mbox{and}\qquad
  \lambda f(x)+(1-\lambda)f(y)>f(t),
}
respectively. Using assertions (ii) and (iv) of \lem{ineq} and that $f$ and $g$ are continuous at $t$, it follows that
\Eq{*}{
  \widetilde{A}^{[2]}_{g}(x,y;\lambda,1-\lambda)
  <t<\widetilde{A}^{[2]}_{f}(x,y;\lambda,1-\lambda).
}
This inequality contradicts \eq{2v}. This contradiction validates the inequality stated in assertion (iv).

Now assume that assertion (iv) holds. Let $t\in C_f\cap C_g$ be fixed and let $x<t<y$ be arbitrary elements of $I$. Observe that the inequality in assertion (iv) can be equivalently written as
\Eq{*}{
  \frac{f(t)-f(x)}{g(t)-g(x)}\geq\frac{f(y)-f(t)}{g(y)-g(t)}.
}
Therefore, for every $t\in C_f\cap C_g$,
\Eq{*}{
  p(t):=\inf_{x\in I,\,x<t}\frac{f(t)-f(x)}{g(t)-g(x)}
  \geq\sup_{y\in I,\,y>t}\frac{f(y)-f(t)}{g(y)-g(t)}.
}
Using the definition of the function $p:C_f\cap C_g\to(0,\infty)$ so defined, for $x,y\in I$ with $x<t<y$, we get
\Eq{*}{
  f(t)-f(x)\geq p(t)(g(t)-g(x)) \qquad\mbox{and}\qquad
  f(y)-f(t)\leq p(t)(g(y)-g(t)).
}
It immediately follows from these inequalities that assertion (v) holds true. 

Next, we verify the implication (v)$\Rightarrow$(i).
Assume assertion (v), and let $n\in\N$ and $x\in I^n$ be fixed. Applying the inequality in assertion (v) with $x:=x_i$ and the adding up the inequalities so obtained, for all $t\in C_f\cap C_g$, we get 
\Eq{*}{
  \sum_{i=1}^n (f(x_i)-f(t))\leq p(t)\sum_{i=1}^n (g(x_i)-g(t)).
}
Choose $t\in C_f\cap C_g$ satisfying the inequality $A^{[n]}_{g}(x)<t$ arbitrarily. Then, by assertion (ii) of \lem{ineq}, we have that
\Eq{*}{
  \frac{g(x_1)+\dots+g(x_n)}{n}<g(t), 
}
whence it follows that $\sum_{i=1}^n (g(x_i)-g(t))<0$. According to the previous inequality, this yields that $\sum_{i=1}^n (f(x_i)-f(t))<0$. Therefore,
\Eq{*}{
  \frac{f(x_1)+\dots+f(x_n)}{n}<f(t). 
}
Using that $f$ is continuous at $t$, by assertion (ii) of \lem{ineq}, we can conclude that $A^{[n]}_{f}(x)<t$. Upon taking the right limit $t\downarrow A_g^{[n]}(x)$ in the dense set $C_f\cap C_g$, we can see that the inequality stated in assertion (i) is valid.

So far, we have established the equivalence of assertions (i)--(v).
Finally, we show that assertions (v), (vi) and (vii) are also equivalent. 

Assume first that (v) holds. For fixed $t\in C_f\cap C_g$, define the function $\varphi_t:\conv(f(I))\to\R$ by
\Eq{*}{
  \varphi_t(u)=\frac{u-f(t)}{p(t)}+g(t) \qquad(u\in\conv(f(I))).
}
Then $\varphi_t$ is a strictly increasing affine function and, by the inequality of assertion (v), we have that
\Eq{xt}{
\varphi_t(f(x))&\leq g(x) \qquad (x\in I,\,t\in C_f\cap C_g)\\
\varphi_t(f(t))&=g(t) \qquad (t\in C_f\cap C_g).
}
Define the extended real-valued function $\varphi:\conv(f(I))\to(-\infty,\infty]$ by
\Eq{*}{
  \varphi(u):=\sup_{t\in C_f\cap C_g}\varphi_t(u) 
  \qquad (u\in\conv(f(I))).
}
Upon taking the supremum in the inequality in \eq{xt}, we obtain that
\Eq{fgx}{
  \varphi(f(x))\leq g(x) \qquad (x\in I).
}
The equality in \eq{xt} implies that
\Eq{*}{
  g(t)=\varphi_t(f(t))
     \leq\varphi(f(t))\leq g(t) \qquad (t\in C_f\cap C_g),
}
which implies that
\Eq{fgt}{
  \varphi(f(t))=g(t) \qquad (t\in C_f\cap C_g).
}
For all $t\in C_f\cap C_g$, we have that $\varphi_t$ is strictly increasing, whence it follows that $\varphi$ is increasing. If $u\in \conv(f(I))$, then there exists $x\in I$ such that $u\leq f(x)$.
Then, using the inequality in \eq{fgx}, we get that $\varphi(u)\leq\varphi(f(x))\leq g(x)$, which shows that $\varphi$ is in fact real-valued. Since $\varphi$ is the supremum of affine functions, it follows that $\varphi$ is convex, which implies that it is also continuous. 

We now are going to point out that $\varphi$ is not only increasing but it is strictly increasing. If this were not be the case, then, using the convexity of $\varphi$, it follows that there exists $u\in\conv(f(I))$ such that $\varphi$ is constant on $J:=\conv(f(I))\cap(-\infty,u]$ and $\varphi$ is strictly increasing on $\conv(f(I))\cap(u,\infty)$. Then, there exists a point $x\in I$ such that $f(x)\leq u$. Let $s,t\in C_f\cap C_g$ with $s<t\le x$. Then, $f(s)<f(t)\leq u$, thus $f(s),f(t)\in J$. Applying that $\varphi$ is constant on $J$, it follows that $\varphi(f(s))=\varphi(f(t))$. On the other hand, in view of the equality \eq{fgt}, the strict increasingness of $g$ implies that
\Eq{*}{
  \varphi(f(s))=g(s)<g(t)=\varphi(f(t)).
}
The contradiction so obtained shows that $\varphi$ must be strictly increasing.

To complete the proof of the implication (v)$\Rightarrow$(vi), we need to show that $\varphi$ is a homeomorphism between the open intervals $\conv(f(I))$ and $\conv(g(I))$. First we establish the inclusion $\varphi(\conv(f(I)))\subseteq\conv(g(I))$. Let $u\in\conv(f(I))$. Then there exist $x,y\in I$ such that 
$f(x)\leq u\leq f(y)$, and there exist $s,t\in C_f\cap C_g$ such that $s\leq x$ and $y\leq t$. These choices, using \eq{fgt}, imply that
\Eq{*}{
  g(s)=\varphi(f(s))\leq \varphi(f(x))\leq\varphi(u)
  \leq \varphi(f(y))\leq \varphi(f(t))=g(t).
}
Therefore,
\Eq{*}{
  \varphi(u)\in[g(s),g(t)]\subseteq\conv(g(I)).
}
To show that the reverse inclusion is also valid, let $v\in\conv(g(I))$. Then there exist $x,y\in I$ such that 
$g(x)\leq v\leq g(y)$, and there exist $s,t\in C_f\cap C_g$ such that $s\leq x$ and $y\leq t$. Now, using \eq{fgt}, we obtain that
\Eq{*}{
  f(s)=\varphi^{-1}(g(s))\leq \varphi^{-1}(g(x))\leq\varphi^{-1}(v)
  \leq \varphi^{-1}(g(y))\leq \varphi^{-1}(g(t))=f(t).
}
Therefore,
\Eq{*}{
  v\in[\varphi(f(s)),\varphi(f(t))]
  =\varphi[f(s),f(t)]\subseteq\varphi(\conv(f(I))).
}

We now prove that assertion (vi) implies (v). Applying that $\varphi$ is convex on the interval $\conv(f(I))$. By well-known properties of convex functions, for every $u\in\conv(f(I))$, we can find a real number $q(u)$ such that
\Eq{*}{
  \varphi(u)+q(u)(v-u)\leq\varphi(v) \qquad (v\in\conv(f(I)))
}
In particular, for all $x\in I$ and $t\in C_f\cap C_g$, with $u:=f(t)$ and $v:=f(x)$, this inequality implies that
\Eq{*}{
  \varphi(f(t))+q(f(t))(f(x)-f(t))\leq\varphi(f(x)) \qquad (x\in I,\,t\in C_f\cap C_g).
}
Using the inequality and the equality stated in (vi) (i.e., the inequality \eq{fgx} and the equality \eq{fgt}), it follows that 
\Eq{*}{
  g(t)+q(f(t))(f(x)-f(t))\leq g(x) \qquad (x\in I,\,t\in C_f\cap C_g).
}
This shows that assertion (v) is validated by the function $p:=\frac{1}{q\circ f}$.

To see that assertions (vi) and (vii) are equivalent, observe that if (vi) holds with an increasing homeomorphism $\varphi$, then (vii) is valid with $\psi:=\varphi^{-1}$ (which is also an increasing homeomorphism). The reverse implication is analogous.
\end{proof}

\section{Equality of generalized quasiarithmetic means}

In this section we show two results concerning equalities of generalized quasiarithmetic means. Both of them are obtained by using a corresponding result for the inequalities proved in the previous section and the easy-to-see fact that equality can be seen as  two combined inequalities.

First result states that, whenever the equality holds between two $n$-variable generalized quasiarithmetic means and $m<n$, then it holds between the corresponding $m$-variable generalized quasiarithmetic means, as well. It is a direct consequence of \thm{M1}.

\Cor{M1}{Let $f,g:I\to\R$ be strictly increasing functions and let $n,m\in\N$ with $m<n$. If, for all $x\in I^n$, the equality
\Eq{*}{
  A_f^{[n]}(x)=A_g^{[n]}(x)
}
holds, then, for all $y\in I^m$,
\Eq{*}{
  A_f^{[m]}(y)=A_g^{[m]}(y).
}}

The following result is a direct consequence of \thm{M2}.

\Cor{M2}{Let $f,g:I\to\R$ be strictly increasing functions. Then the following assertions are equivalent.
\begin{enumerate}[(i)]
 \item For all $n\in\N$ and $x\in I^n$,
 \Eq{*}{
   A^{[n]}_{f}(x)=A^{[n]}_{g}(x).
 }
 \item For all $n\in\N$, $x\in I^n$, and $\lambda\in\Lambda_n$
 \Eq{*}{
   \widetilde{A}^{[n]}_{f}(x;\lambda)=\widetilde{A}^{[n]}_{g}(x;\lambda).
 }
 \item There exists $n\in\N\setminus\{1\}$ such that, for all $x\in I^n$ and $\lambda\in\Lambda_n$,
 \Eq{*}{
   \widetilde{A}^{[n]}_{f}(x;\lambda)=\widetilde{A}^{[n]}_{g}(x;\lambda).
 }
 \item For all $x,y$ in $I$ with $x<y$ and $t\in (x,y)\cap C_f\cap C_g$,
 \Eq{*}{
  \frac{f(y)-f(t)}{f(y)-f(x)}=\frac{g(y)-g(t)}{g(y)-g(x)}.
 }
 \item There exist two constants $\alpha>0$ and $\beta$ such that, for all $x\in I$,
 \Eq{*}{
   f(x)=\alpha g(x)+\beta.
 }
\end{enumerate}}

\begin{proof}
The equivalence of assertions (i)--(iv) directly follows from the equivalence of assertions (i)--(iv) of \thm{M2}.

Assume now that assertion (i) is valid. Then, applying that assertion (i) implies assertion (v) in \thm{M2} twice, we get that there exist two positive functions $p,q:(C_f\cap C_g)\to\R_+$ such that, for all $x\in I$ and $t\in C_f\cap C_g$,
\Eq{2in}{
   f(x)-f(t)\leq p(t)(g(x)-g(t)) \qquad\mbox{and}\qquad
   g(x)-g(t)\leq q(t)(f(x)-f(t)).
}
Therefore, for all $x\in I$ and $t\in C_f\cap C_g$, we have that
\Eq{*}{
  \frac{1}{q(t)}(g(x)-g(t))\leq p(t)(g(x)-g(t)).
}
Choosing first $x$ to be smaller than $t$ and then to be bigger than $t$, it follows that $\frac{1}{q(t)}=p(t)$ for all $t\in C_f\cap C_g$. Thus, the two inequalities in \eq{2in} imply, for all $x\in I$ and $t\in C_f\cap C_g$, that
\Eq{*}{
  f(x)-f(t)=p(t)(g(x)-g(t))
}
Let $t\in C_f\cap C_g$ be a fixed element and define $\alpha:=p(t)$ and $\beta:=f(t)-p(t)g(t)$. Then the above equality shows that assertion (v) is valid.

Conversely, assume that assertion (v) holds for some real constants $\alpha>0$ and $\beta$. Then assertion (iv) follows immediately.
\end{proof}

\section{An example and concluding remarks}

Define, for $n\in\N$, the mean $\mathfrak{M}^{[n]}:\R^n\to\R$ as follows:
\Eq{*}{
  \mathfrak{M}^{[n]}(x)
  =\begin{cases}
     0 &\mbox{if } \min(x)\leq 0\leq\max(x),\\[2mm]
     \dfrac{x_1+\dots+x_n}{n} &\mbox{otherwise}.
    \end{cases}
}
Similar constructions of means were introduced in the papers \cite{DarJarJar16} and \cite{JarJar17}. In what follows, we investigate the representability of the mean $\mathfrak{M}^{[n]}$ as a (generalized) quasiarithmetic mean on the subinterval $I$. It is obvious that if $0\not\in I$, then $\mathfrak{M}^{[n]}$ restricted to $I^n$ is equal to the $n$-variable arithmetic mean. Therefore, the interesting case is when $0\in I$.

\Prp{M}{Assume that $0\in I$. Then, we have the following assertions:
\begin{enumerate}[(i)]
 \item For all $n\in\N\setminus\{1\}$, the mean $\mathfrak{M}^{[n]}$ restricted to $I^n$ is not (standard) quasiarithmetic mean on $I^n$, i.e., there is no strictly increasing and continuous function $f:I\to\R$ such that 
 \Eq{M}{
   \mathfrak{M}^{[n]}(x)=A_f^{[n]}(x)\qquad(x\in I^n)
 }
 is valid.
 \item For all $n\in\N\setminus\{1\}$, the mean $\mathfrak{M}^{[n]}$ restricted to $I^n$ is a generalized quasiarithmetic mean, that is, there exists a strictly increasing function $f:I\to\R$ such that \eq{M} holds, if and only if $I$ is bounded.
 \item There is no strictly increasing function $f:I\to\R$ such that, for all $n\in\N$, the equality \eq{M} holds.
\end{enumerate}
}

\begin{proof} 
 
To show the first assertion it is sufficient to observe that all (standard) quasiarithmetic means are strictly increasing (as stated in assertion (iii) of \thm{C}) while the mean $\mathfrak{M}^{[n]}$ does not possess this property if $n\geq2$.

Now we proceed to the proof of assertion (ii). 
Let $n\in\N\setminus\{1\}$ be fixed. First assume that $I$ is bounded and let $a<0<b$ denote the endpoints of $I$. Define $f:[a,b]\to\R$ by
\Eq{f}{
  f(x):=\begin{cases}
         -\frac{x}{an}-1 &\mbox{if } x\in[a,0)\\
         0 &\mbox{if } x=0,\\
         \frac{x}{bn}+1 &\mbox{if } x\in(0,b]\\
        \end{cases}\qquad(x\in I).
}
It is easy to see that $f$ is strictly increasing. 
For the generalized inverse of $f$ (restricted to $(a,b)$), we obtain
\Eq{if}{
  f^{(-1)}(u)
  =\begin{cases}
    -an(u+1) & \mbox{if } u\in\big(-\frac{n+1}{n},-1\big),\\
    0 & \mbox{if } u\in[-1,1],\\
    bn(u-1) & \mbox{if } u\in\big(1,\frac{n+1}{n}\big).
   \end{cases}
}
We are going to show that \eq{M} holds for all $x\in I^n$. Let $x\in I^n$ be fixed.

Observe that $f$ is affine on the open subintervals $(a,0)$ and $(0,b)$ and hence, in these subintervals $f$ generates the arithmetic mean. Therefore, if either $\max(x)<0$ or $0<\min(x)$, then \eq{M} is valid. 

It remains to show that \eq{M} is satisfied in the case $\min(x)\leq 0\leq\max(x)$, that is, $A_f^{[n]}(x)=0$. Equivalently, in view of formula \eq{if}, we have to prove that
\Eq{R}{
  \frac{f(x_1)+\dots+f(x_n)}{n}\in[-1,1].
}
Due to the inequality $\min(x)\leq0$, for some $i$ we have that $x_i\leq0$, and $x_j\leq b$ for all $j\neq i$. Therefore,
\Eq{*}{
  \frac{f(x_1)+\dots+f(x_n)}{n}
  \leq \frac{f(0)+(n-1)f(b)}{n}
  =\frac{(n-1)\frac{n+1}{n}}{n}
  =\frac{n^2-1}{n^2}<1.
}
A similar argument shows that
\Eq{*}{
  \frac{f(x_1)+\dots+f(x_n)}{n}
  \geq \frac{f(0)+(n-1)f(a)}{n}
  =\frac{-(n-1)\frac{n+1}{n}}{n}
  =\frac{1-n^2}{n^2}>-1.
}
Thus we have proved that \eq{R} is valid and hence $A_f^{[n]}(x)=0$ holds, which completes the proof of \eq{M} in the case $\min(x)\leq 0\leq\max(x)$. Consequently, $\mathfrak{M}^{[n]}$ restricted to $I^n$ is a generalized quasiarithmetic mean provided that $I$ is bounded.

To complete the proof of assertion (ii) assume now that $I$ is unbounded. To the contrary of the statement, assume that there exists a strictly increasing function $f:I\to\R$ such that \eq{M} is valid for some $n\geq2$. Consider now the subcase when $\sup I=\infty$.
Then, $A_f^{[n]}$ restricted to the interval $(0,\infty)^n$ is equal to the arithmetic mean, which, as a quasiarithmetic mean is generated by an affine function. Therefore, according to assertion (v) of \cor{M2} below, there exist two real constants $\alpha>0$ and $\beta$ such that $f(t)=\alpha t+\beta$ for all $t>0$. Then $\beta=f_+(0)\geq f(0)$, and one can get that $f^{(-1)}(u)=\frac{1}{\alpha}(u-\beta)>0$ holds for all $u>\beta$.

On one hand, for arbitrary $t>0$, we have that 
\Eq{*}{
  \mathfrak{M}^{[n]}(0,\dots,0,t)=0.
} 
On the other hand, for $t>\frac{(n-1)(\beta-f(0))}{\alpha}\geq0$, we have that $\frac{(n-1)f(0)+\alpha t+\beta}{n}>\beta$. Therefore
\Eq{*}{
  A_f^{[n]}(0,\dots,0,t)
  =f^{(-1)}\bigg(\frac{(n-1)f(0)+f(t)}{n}\bigg)
  =f^{(-1)}\bigg(\frac{(n-1)f(0)+\alpha t+\beta}{n}\bigg)
  >0,
}
which shows that \eq{M} cannot hold at $x=(0,\dots,0,t)$ with $t>\frac{(n-1)(\beta-f(0))}{\alpha}$.

To verify the last assertion, assume that \eq{M} holds for all $n\in\N$ with some strictly increasing function $f:I\to\R$. Then, according to assertion (iii) of \prp{QAM2}, for all $x,y\in I$, we have that 
\Eq{*}{
  \lim_{n \to \infty} A_f^{[n]}\big(x,\underbrace{y,\dots,y}_{(n-1)\text{ times}}\big)
  =\lim_{n \to \infty} \widetilde{A}_f^{[2]}\big(x,y;\tfrac1n,1-\tfrac1n\big)
  =\widetilde{A}_f^{[2]}\big(x,y;0,1\big)
  =y,
}
while, for all $x<0<y$ and $n\geq2$, we have that
\Eq{*}{
\mathfrak{M}^{[n]}\big(x,\underbrace{y,\dots,y}_{(n-1)\text{ times}}\big)=0.
}
This equality causes a contradiction in \eq{M} for large $n$. 
\end{proof}

\subsection*{Concluding remarks} It this section we would like to emphasize few important remarks concerning generalized quasiarithmetic means (which we will call here briefy g.q.a. means):
\begin{enumerate}[(1)]
 \item G.q.a. means are the generalizations of quasiarithmetic means, however they also form a subfamily of semideviation means introduced by P\'ales \cite{Pal89b} and developed in P\'ales--Pasteczka \cite{PalPas19a}. In~particular, all results concerning this family are also applicable to g.q.a.\ means.
 \item Each Kolmogorov's axiom expressed in \thm{C} is either valid for g.q.a.\ mean (reflexivity and symmetry) or they imply that a g.q.a.\ mean is a quasiarithmetic mean (such as continuity, strict increasingness, associativity).
 \item An $n$-variable g.q.a.\ mean completely determines the corresponding $m$-variable mean whenever $m<n$, however the converse of this statement is not valid.
\item Every g.q.a.\ mean possesses a lower semicontinuous and upper semicontinuous envelop, which are also g.q.a.\ means.
\item Except a few results, weighted g.q.a.\ means are (intentionally) a bit aside from the main plot of this paper. Since the construction of weighted g.q.a.\ means follows the general principle introduced in \cite{PalPas18b}, this family can be reconstructed using the nonweighted g.q.a.\ means which imply their properties in a natural manner.
\item Due to \thm{M2} part (vi), we know that the comparability (for vectors of arbitrary lengths) of two g.q.a.\ means is localizable (i.e.\ they are comparable if, and only if, they are comparable for arguments from a neighbourhood of the diagonal). On the other hand, this is not the case when we restrict to vectors with a fixed number of elements.
\item The family of g.q.a.\ means does not follow any extension principles. More precisely:
\begin{itemize}
 \item[$\circ$] For every $n\in\N\setminus\{1\}$, the mean $\mathfrak{M}^{[n]} \colon \R^n \to \R$ constructed in this section is not a g.q.a.\ mean, however its restriction $\mathfrak{M}^{[n]}|_{I^n}$ is a g.q.a.\ mean for every compact subinterval $I \subset \R$;
 \item[$\circ$] The mean $\mathfrak{M} \colon \bigcup_{n=1}^\infty I^n \to I$ is not a g.q.a.\ mean, however its restriction $\mathfrak{M}|_{I^n}$ is a g.q.a.\ mean for every $n \in \N$.
\end{itemize}
\item The family of g.q.a.\ means can be extended to the integral setting in the spirit of Hardy--Littlewood--P\'olya \cite{HarLitPol34}. Namely, for a strictly monotone function $f \colon I \to \R$ and a function $\varphi \colon J \to I$, we can define
\Eq{*}{
\mathscr{A}_f(\varphi):=f^{(-1)} \Big(\frac{1}{|J|}\int_J f\circ\varphi(x)\:dx\Big).
}
At this stage the mutual relationship between such means and integral quasiarithmetic means is unknown.
\end{enumerate}

\subsection*{Conflict of interest} On behalf of all authors, the corresponding author states that there is no conflict of interest.

\subsection*{Data availability statement}  Not applicable.


\end{document}